\documentclass{article}
\usepackage{amsmath,amssymb}
  \usepackage[noend]{algorithmic}
  \usepackage{cite}

  \newtheorem{lemma}{Lemma}
  \newtheorem{theorem}{Theorem}
  
  \newtheorem{asm}{Assumption}
  \newtheorem{alg}{Algorithm}

  \def\Rset{\mathbb{R}}
  \def\Nset{\mathbb{N}}
  \def\proof{\noindent\mbox{\bf{Proof}\ }}
  \def\qed{\ \hfill $\square$}
  \newcommand{\norm}[1]{\| #1 \|}
  
  \def\rank{\mbox{rank\,}}
  \def\algbreak{\textbf{break}}

\begin{document}

  \title{Adaptive set-point regulation of discrete-time nonlinear systems
    \thanks{
      This work was supported by the Japan Society for the Promotion of 
      Science under Grant-in-Aid for Scientific Research (C) 23560535.
      This manuscript is a former version of the manuscript the author has submitted to 
      International Journal of Adaptive Control and Signal Processing.  
      The manuscript was rejected, and a revision is in preparation, but 
      this version does not reflect the comments of the referees to the rejected version.  
    }
  }
  \author{Shigeru~Hanba
    \thanks{
      Department of Electrical and Electronics Engineering, University of the Ryukyus,
      1 Senbaru Nishihara, Nakagami-gun, Okinawa 903-0213, Japan; email: sh@gargoyle.eee.u-ryukyu.ac.jp
    }
  }
  \maketitle

  \begin{abstract}
  In this paper, adaptive set-point regulation controllers for discrete-time nonlinear systems are constructed.  
  The system to be controlled is assumed to have a parametric uncertainty, 
  and an excitation signal is used in order to obtain the parameter estimate.  
  The proposed controller belongs to the category of indirect adaptive controllers, 
  and its construction is based on the policy of calculating the control input rather than that of obtaining a control law.  
  The proposed method solves the adaptive set-point regulation problem
  under the (possibly minimal) assumption that the target state is reachable provided that the parameter is known.  
  Additional feature of the proposed method is that Lyapunov-like functions have not been used in the construction of the controllers.  
  \end{abstract}

  \paragraph*{keywords}
  adaptive set-point regulation, discrete-time nonlinear systems, persistent excitation

 \section{Introduction}
 For decades, 
 adaptive control of nonlinear systems has been an active area of research,
 and several design methods have been established
 for both continuous-time and discrete-time systems
 \cite{Seto1994,Krstic1995,Xie1999,Lin2002,Loh2003,Tyukin2004,Hayakawa2004,Zhang2007,Postoyan2007,Ge2008,
   Hartwich2008,Wu2011,Yalcin2012,Dai2013,Li2014,Tao2014}.  
 A majority of these design methods first assume that the systems are described in some canonical forms
 and that parametric Lyapunov-like functions are known; 
 they then construct controllers together with tuners of specific forms 
 to obtain sufficient conditions for the stability of the closed-loop system.  
 In other words, their sufficient conditions are ``method driven. ''
 Therefore the question naturally arises: 
 what is a nearly minimal sufficient condition for a nonlinear system 
 with parametric uncertainty to permit stable adaptive controllers?  
 The present paper is an attempt to answer this question.  
 
 The objective of this paper is to construct adaptive controllers for discrete-time nonlinear systems 
 that drive the state of the system into a neighborhood of a ``target state'' by finite-time control
 under the assumption that the target state is finite-time reachable if the parameter is known.  
 We assume a certain kind of parameter identifiability (the precise statement is given below)
 together with an excitation signal, and we construct the controllers by 
 following the policy of calculating the control rather than that of obtaining a control law, 
 assuming that numerical solutions of nonlinear equations are available either exactly or with the desired accuracy.  
 The proposed method is not supposed to be used in consecutive operations 
 --- it is assumed that the controllers terminate if the state reaches
 a neighborhood of the target state, and the parameter estimate is used for other purpose
 (e.g. for parametric local stabilizing controller).  

 \section{Definitions and notations}
 Consider a discrete-time nonlinear system with parametric uncertainty of the form 
 \begin{equation}
  x(t+1)=f(x(t),u(t),\theta), 
  \label{eq:system}
 \end{equation}
 where $x(t) \in \Rset^{n}$ is the state, 
 $u(t) \in \Rset^{n_u}$ is the control input, 
 and $\theta \in \Rset^{n_\theta}$ is the parameter to be estimated.  
 The parameter $\theta$ is assumed to be inside a compact and convex set $\Omega_{\theta}$, 
 and the function $f$ is assumed to be $C^1$ with respect to all arguments.  

 Henceforth, we use the following notations. 
 We denote the sequence of inputs $(u(t_0),\ldots,u(t_1))$ by $u[t_0,t_1]$.  
 The sequence of the state 
 $(x(t_0),\ldots,x(t_1))$ is denoted as $X[t_0,t_1;u,\theta]$,
 where the symbols $u,\theta$ have been added to emphasize their effect.  
 Although $u[t_0,t_1]$ and $X[t_0,t_1;u,\theta]$ are sequences of vectors with length $t_1-t_0+1$, 
 we sometimes identify them with vectors in $\Rset^{n_u(t_1-t_0+1)}$ and $\Rset^{n(t_1-t_0+1)}$,
 respectively.  

 The solution of (\ref{eq:system}) at $t=t_1$
 initialized at $t=t_0$ with $x(t_0)$ is denoted by $\varphi(t_1,t_0,x(t_0);u,\theta)$.  
 The symbols $B(x,\rho)$ and $\overline{B}(x,\rho)$ denote the open
 and closed balls centered at $x$ with radius $\rho$.
 For a sequence $u[t_0,t_1]$,
 $\overline{B}^{\infty}(u[t_0,t_1],\rho)$ denotes the set
 $\{v[t_0,t_1]: \forall t \in \{t_0,\ldots,t_t\}, v(t) \in \overline{B}(u(t),\rho)\}$.
 The symbol $\Nset$ denotes the set of nonnegative integers.  

 Let the target state be $x_{\ast}$.  
 Henceforth, we assume that the target state is finite-time reachable
 in the following sense, which is a parametric counterpart of 
 those given in \cite{Hanba2009a}.  
 \begin{asm}\label{asm:reachability}
   $\forall x$, $\forall \theta \in \Omega_{\theta}$, 
   $\exists N>0$, $\exists u[0,N-1]$, 
   \begin{enumerate}
   \item $\varphi(N,0,x;u,\theta)=x_{\ast}$, 
   \item 
     $\rank \frac{\partial \varphi(N,0,x;u,\theta)}{\partial u[0,N-1]}=n$.  
   \end{enumerate}
 \end{asm}
 The first expression in Assumption~\ref{asm:reachability} is the algebraic reachability, 
 and the second expression is the nonlinear counterpart of the controllability rank condition.  
 Combining them implies a certain kind of 
 ``uniform controllability''\cite{Hanba2009a}, 
 which has an observability counterpart\cite{Hanba2009b,Hanba2010}.  

 As for the parameter identification, 
 we assume the existence of the following excitation signal.  
 \begin{asm}\label{asm:excitation}\rm
  $\exists N$, $\exists u[0,N-1]$, 
  $\forall x$, $\forall \theta \in \Omega_{\theta}$, 
  $\rank \frac{\partial X[0,N;u,\theta]}{ \partial \theta}=n_{\theta}$. 
 \end{asm}

 In what follows, 
 we construct a controller with structure similar to block model predictive controllers\cite{Hanba2009a,Sun2007a}
 together with parameter estimators based on nonlinear equation solvers. 
 Although model predictive control is an application-oriented method, 
 the use of this `block model predictive control' structure in our adaptive controller
 has completely different objective.  It is used as a theoretical tool to show that
 Assumption~\ref{asm:excitation} serves as a `persistent excitation condition' 
 for a nonlinear system of the form (\ref{eq:system}).  
 Basically, our control strategy is as follows.  
 Partition the time interval $\Nset$ into blocks of finite length 
 (the length of each block is determined adaptively. )
 Let $t=T_k$ be the beginning of the $k$-th block.  
 At this time instant, 
 update the parameter estimate $\theta(T_k)$ using the entire sequence of past states.  
 Then, obtain $N_k$ and $u[T_k,T_k+N_k-1]$ such that 
 $\varphi(T_k+N_k,T_k,x(T_k);u,\theta(T_k))=x_{\ast}$
 or  $\varphi(T_k+N_k,T_k,x(T_k);u,\theta(T_k))\in {\cal N}(x_{\ast})$, 
 where ${\cal N}(x_{\ast})$ is a neighborhood of $x_{\ast}$.  
 Next, apply the input sequence $u[T_k,T_k+N_k-1]$ to the system (\ref{eq:system})
 until $t=T_k+N_k-1$ is in an open-loop fashion.  

 To avoid the overuse of subscripts, 
 henceforth, 
 we employ the following simplified notations.  
 First, $\theta(T_k)$ is rewritten as $\theta_k$.  
 For the parameter estimation, the entire sequence of states up to $T_k$, $X[1,T_k;u,\theta]$ will be used, 
 but what really matters is the dependence on $\theta$ only; 
 we rewrite this expression in the column vector form and let $g_k(\theta)=( (x(1))^T,\ldots,(x(T_k))^T)^T$, 
 omitting unnecessary variables to avoid confusion.  
 Similarly, we rewrite $\varphi (T_k+N_k,0,x(0);u,\theta)$ as $h_k(\theta)$.  

 \section{Main results}

 We first consider the ideal case where solutions of nonlinear equations 
 are available exactly, and we then consider the case where numerical errors 
 to the solutions of nonlinear equations do exist.  
 It is to be emphasized, however, that 
 the assumption that an exact solution of the nonlinear equation of the parameter estimate
 is available, that is, $g_k(\theta_k)=g_k(\theta)$, does not always imply that $\theta_k=\theta$
 because $g_k(\theta)$ is not always a global injection ---
 Assumption~\ref{asm:excitation} merely assures that it is a local injection.  

 The first algorithm of adaptive set-point regulation is as follows.  
 \begin{alg}\label{alg:excitation_zero_error}\rm
   \ \\
   \noindent{\bf{(Initialization)}}
   Given $x(0)$, let $u[0,N_0-1]$ be the excitation signal,
   $T_0=0$, $T_1=N_0$, $k=1$; 
   apply $u[0,N_0-1]$to system (\ref{eq:system}) to obtain $x(T_1)$.  

   \noindent{\bf{(Loop)}}
   \begin{algorithmic}
     \IF {$x(T_k)=x_{\ast}$}
     \STATE \algbreak;
     \ELSE
     \STATE Obtain $\theta_k \in \Omega_{\theta}$ that satisfy $g_k(\theta_k)=g_k(\theta)$;
     \STATE Obtain $u[T_k,T_k+N_k-1]$ that satisfy $h_k(\theta_k) =x_{\ast}$;
     \STATE Let $T_{k+1}=T_{k}+N_k$;
     \STATE Apply $u[T_k,T_{k+1}-1]$ to the system (\ref{eq:system}) to obtain $x(T_{k+1})$.  
     \STATE Let $k=k+1$;
     \ENDIF
   \end{algorithmic}
 \end{alg}

 \begin{theorem}\label{thm:zero_error}
   Under the first condition of Assumption~\ref{asm:reachability}
   and Assumption~\ref{asm:excitation}, 
   Algorithm~\ref{alg:excitation_zero_error}
   terminates after finitely many iterations, 
   and the state of (\ref{eq:system}) reaches to the target state $x_{\ast}$
   provided that exact solutions of nonlinear equations are available.  
 \end{theorem}

 \proof
 First, note that our assumptions permit that 
 Algorithm~\ref{alg:excitation_zero_error} is always feasible.  

 We prove our assertion by contradiction.  
 Suppose that Algorithm~\ref{alg:excitation_zero_error} never terminates after finitely many iterations.  
 Then, the resulting sequence of the parameter estimate $(\theta_k)_{k \in \Nset}$ 
 is an infinite sequence in the compact set $\Omega_{\theta}$ and hence,
 it has at least one limit point.  Let $\theta_{\sharp}$ be one of 
 its limit points.  

 We first prove that $\forall j, g_j(\theta_{\sharp})=g_j(\theta)$.  
 Because we have assumed that exact solutions to nonlinear equations are available, 
 $\forall j$, $g_j(\theta_j)=g_j(\theta)$.  
 Moreover, for $k_1 <k_2$, the relation between  
 $g_{k_1}(\theta)$
 and 
 $g_{k_2}(\theta)$
 are given by 
 \[
 g_{k_2}(\theta)=
 \begin{pmatrix}
   g_{k_1}(\theta)\\
   x(T_{k_1}+1)\\
   \vdots\\
   x(T_{k_2})
 \end{pmatrix}.  
 \]
 Due to this structure, 
 we call that $g_{k_1}(\theta)$ is an initial segment of $g_{k_2}(\theta)$.  
 Because $\theta_{\sharp}$ is a limit point, 
 there is a subsequence 
 $(\theta_{k_l})_{l \in \Nset}$
 of 
 $(\theta_k)_{k \in \Nset}$
 that converges to $\theta_{\sharp}$.  
 For any $j >0$, $g_j(\theta)$ is a continuous function of $\theta$, 
 and $\forall l$ such that $j \leq k_l$,  
 $g_j(\theta_{k_l})=g_j(\theta)$ 
 because  $g_j(\theta_{k_l})$ is an initial segment of $g_{k_l}(\theta)$.  
 Since $\lim_{l\rightarrow\infty} \theta_{k_l}=\theta_{\sharp}$
 and $g_j$ is continuous,  $g_j(\theta_{\sharp})=g_j(\theta)$.  

 Next, we prove by contradiction that $\exists \overline{l}$, $\forall l \geq \overline{l}$, $\theta_{k_l} = \theta_{\sharp}$.  
 Suppose that $\forall \overline{l}$, $\exists l \geq \overline{l}$, $\theta_{k_l} \neq \theta_{\sharp}$.  
 Then, $\forall \varepsilon>0$,  $\exists \theta_{k_l}$, $\theta_{k_l} \in B(\theta_{\sharp},\varepsilon) \setminus \{\theta_{\sharp}\}$.  
 Because $g_1$ (the first segment of $g_{k_l}$) is $C^1$ and $\frac{\partial g_1}{\partial \theta}$ is of full rank, 
 \begin{equation}
   \exists c>0, \exists \varepsilon >0, 
   \left ( \norm{\theta_{k_l}-\theta_{\sharp}}<\varepsilon \Rightarrow 
   \norm{g_1(\theta_{k_l})-g_1(\theta_{\sharp})} \geq c \norm{\theta_{k_l}-\theta_{\sharp}}\right ).  
   \label{eq:Jacobian_bound}
 \end{equation}
 This contradicts the assumption that $g_1(\theta_{k_l})=g_1(\theta)=g_1(\theta_{\sharp})$.  
 Hence, $(\theta_{k_l})_{l \in \Nset}$ converges to $\theta_{\sharp}$ after finitely many iterations,
 and $\forall l \geq \overline{l}$, $\theta_{k_l}=\theta_{\sharp}$.  This also implies that
 there are infinitely many $k$ such that $\theta_{k}=\theta_{\sharp}$.  

 Let $k$, $k^{\prime}$ be such that $k < k^{\prime}$ and $\theta_k=\theta_{k^{\prime}}=\theta_{\sharp}$.  
 We conclude our analysis by showing that $x(T_{k+1})=x_{\ast}$.  
 To see this, we recall our parameter tuning and control mechanism.  
 At the beginning of the $k$-th block, 
 the parameter estimate is updated from $\theta_{k-1}$ to $\theta_{k}$
 to satisfy $g_k(\theta_{k})=g_k(\theta_{\sharp})=g_k(\theta)$.  
 The predicted trajectory based on $\theta_{k}$ is $g_{k+1}(\theta_{k})$, 
 and the control input is determined to make the last $n$ components of $g_{k+1}(\theta_{k})$, 
 $\varphi(T_{k+1},T_k,x(T_k);u,\theta_k)$, 
 identical to the target state $x_{\ast}$.  
 At this state, superficially, 
 it is not assured that $g_{k+1}(\theta_{k})=g_{k+1}(\theta)$, where 
 $g_{k+1}(\theta)$ corresponds to the actual trajectory.  
 However, 
 $g_{k^{\prime}}(\theta_{k^{\prime}})= g_{k^{\prime}}(\theta_{\sharp})=g_{k^{\prime}}(\theta)$, 
 and 
 $g_{k+1}(\theta_{k^{\prime}})= g_{k+1}(\theta_{\sharp})=g_{k+1}(\theta)$
 is its initial segment.  Therefore, 
 $\varphi(T_{k+1},T_k,x(T_k);u,\theta_k)=\varphi(T_{k+1},T_k,x(T_k);u,\theta)=x_{\ast}$.  
 This is a contradiction because 
 we have supposed that
 the algorithm
 never terminates after finitely many iterations. \qed

 Next, we consider the case where solutions to nonlinear equations
 may contain numerical errors, that is, 
 a numerical solution to a nonlinear equation $g(\theta)=0$
 (we temporally denote it by $\widehat{\theta}$) satisfies $\norm{g(\widehat{\theta})-g(\theta)}\leq \varepsilon$
 for some $\varepsilon>0$, but the size of $\varepsilon$ may be arbitrarily specified 
 by a numerical nonlinear equation solver --- generally, such specification is possible
 by adequately tuning the termination condition of the solver, as far as the CPU time permits it.  

 In Algorithm~\ref{alg:excitation_zero_error}, where we have assumed exact solutions to nonlinear equations, 
 there has been no limitation on the length of the blocks and the amplitude of the inputs.  
 They may be arbitrary, and the ``exact solution'' assumption absorbs all of their effect.  
 In contrast, for inexact solution cases, they should be upper-bounded by some constant.  
 The existence of the upper bound (and hence, feasibility) is assured by the following lemma, 
 which is a variant of Lemma~2 in \cite{Hanba2009a}.  

 \begin{lemma}
   Under Assumption~\ref{asm:reachability}, 
   for a fixed $x$, 
   the length of the control block and the amplitude of the control inputs
   that drive $x$ into the target state $x_{\ast}$ 
   are uniformly bounded for all admissible parameters in $\Omega_\theta$
   in the following sense: 
   $\exists N_x >0$, $\exists \rho_x>0$, 
   $\forall \theta \in \Omega_{\theta}$, 
   $\exists N \leq N_x$, 
   $\exists u[0,N-1] \in \overline{B}^{\infty}(0,\rho_x)$, 
   $\varphi(N,0,x;u,\theta)=x_{\ast}$.  
 \end{lemma}

 \proof  
 The proof is similar to that of Lemma~2 of \cite{Hanba2009a}
 and hence, it is omitted.  
 \qed

 Henceforth, we assume the following.  
 \begin{asm}\label{asm:N_x and rho_x}
   For each $x$, 
   $N_x$ and $\rho_x$ are known a priori.     
 \end{asm}

 Our algorithm based on inexact numerical solution also applies
 the excitation signal of Assumption~\ref{asm:excitation} to the system (\ref{eq:system})
 at the beginning of the first control block.  We have not yet described the algorithm itself, 
 but the function $g_1(\theta)$ of Algorithm~\ref{alg:excitation_zero_error} 
 is independent of the algorithm and hence is already determined.  
 In the proof of Theorem~\ref{thm:zero_error}, we have used the fact that 
 for a fixed $\theta_{\sharp}$, (\ref{eq:Jacobian_bound}) holds 
 because $\rank \frac{\partial g_1}{\partial \theta} = n_{\theta}$.  
 Our inexact numerical solution counterpart requires 
 its ``uniform counterpart. ''
 \begin{lemma}\label{lem:uniform_jacobian}
   $\exists \varepsilon_{g_1}>0$, $\exists c_{g_1}>0$,
   $\forall \theta_1,\theta_2 \in \Omega_{\theta}$,
   \[
   \norm{\theta_1-\theta_2} \leq \varepsilon_{g_1} \Rightarrow  
   \norm{g_1(\theta_1)-g_1(\theta_2)} \geq c_{g_1} \norm{\theta_1-\theta_2}.     
   \]
 \end{lemma}

 \proof
 The proof is by contradiction.  
 Suppose that $\forall \varepsilon>0$, $\forall c>0$,
 $\exists \theta_1, \theta_2 \in \Omega_{\theta}$,
 $\norm{\theta_1-\theta_2} < \varepsilon$ and 
 \begin{equation}
   \norm{g_1(\theta_1)-g_1(\theta_2)} < c \norm{\theta_1 - \theta_2}.  
   \label{eq:dx_lb_001}
 \end{equation}

 Let 
 $Jg_1=\frac{\partial g_1}{\partial \theta}$
 and $\lambda_{\rm{min}}=\min \{\norm{(Jg_1)(\theta)v}: v \in \Rset^{n_{\theta}}, \norm{v}=1; \theta \in \Omega_{\theta} \}$.  
 Because $Jg_1$ is continuous and of full rank, $\lambda_{\rm{min}}>0$.  
 Let $c=\frac{\lambda_{\rm{min}}}{2}$, 
 and let $(\theta_1(k),\theta_2(k))$ be the pair in 
 $\Omega_{\theta}$ that satisfies (\ref{eq:dx_lb_001}) for $\varepsilon=1/k$.  
 Because (\ref{eq:dx_lb_001}) does not include equality,
 $\theta_1(k) \neq \theta_2(k)$.  
 Because $\Omega_{\theta} \times \Omega_{\theta}$ is compact, 
 $(\theta_1(k),\theta_2(k))_{k \in \Nset}$ has an accumulation point  $(\theta_{\sharp},\theta_{\sharp})$.  
 By Taylor's formula and the assumption that $\Omega_{\theta}$ is convex, 
 $g_1(\theta_2)-g_1(\theta_1)=(Rg_1)(\theta_1,\theta_2,p)(\theta_2-\theta_1)$, 
 where $p =(p_1,\ldots,p_{n_{\theta}}) \in \prod^{n_{\theta}}[0,1]$ and 
 \[
 (Rg_1)(\theta_1,\theta_2,p)
 =
 \begin{pmatrix}
  \left . \frac{\partial g_{1,1}}{\partial \theta} \right |_{p_1 \theta_2+(1-p_1)\theta_1}\\
  \ldots\\
  \left . \frac{\partial g_{1,{n_{\theta}}}}{\partial \theta}  \right |_{p_{n_\theta} \theta_2+(1-p_{n_{\theta}})\theta_1}
 \end{pmatrix}.  
 \]
 Because $Rg_1$ is continuous and its domain is compact, 
 it is uniformly continuous, 
 and $(Rg_1)(\theta_{\sharp},\theta_{\sharp},p)=(Jg_1)(\theta_{\sharp})$.  
 Therefore, $\forall \varepsilon>0$, $\exists \delta>0$, 
 \begin{equation}
   \max_{i=1,2}\{\norm{\theta_i-\theta_{\sharp}}\}<\delta \Rightarrow
   \norm{(Rg_1)(\theta_1,\theta_2,p)-(Jg_1)(\theta_{\sharp})}<\varepsilon. 
   \label{eq:uniform_jacobian_c}
 \end{equation}
 Let $(\theta_1(k),\theta_2(k))$ be the pair that satisfies
 (\ref{eq:uniform_jacobian_c}) for $\varepsilon<\frac{\lambda_{\rm{min}}}{2}$.  
 Then, since 
 \[
 \begin{split}
 g_1(\theta_2(k))-g_1(\theta_1(k))
 &=(Jg_1)(\theta_{\sharp})(\theta_2(k)-\theta_1(k))\\
 &+((Rg_1)(\theta_1,\theta_2,p)-(Jg_1)(\theta_{\sharp}))(\theta_2(k)-\theta_1(k)), 
 \end{split}
 \]
 it follows that $\norm{g_1(\theta_2(k))-g_1(\theta_1(k))} \geq \frac{\lambda_{\rm{min}}}{2} \norm{\theta_2(k)-\theta_1(k)}$, 
 contradicting (\ref{eq:dx_lb_001}).  
 \qed

 Now, we describe the algorithm.  
 In our algorithm, 
 the numerical error of the solutions of nonlinear equations are
 treated by a method that is similar to the trust-region method of nonlinear programming\cite{Conn2000}.  

 \begin{alg}\label{alg:excitation_finite_error}\rm
   \ \\
  \noindent{\bf{(Initialization)}}
  Given $x(0)$, 
  choose a constant $\beta$ $(0<\beta<1)$, 
  $\mu_0>0$, $\kappa_0>0$, and $\varepsilon_{\rm{fin}} >0$.  
  Let $u[0,N_0-1]$ be the excitation signal,
  $T_0=0$, $T_1=N_0$, $k=1$; 
  Apply $u[0,N_0-1]$ to the system (\ref{eq:system}) to obtain $x(T_1)$.  

  \noindent{\bf{(Loop)}}
  \begin{algorithmic}
    \IF {$\norm{x(T_k)} < \varepsilon_{\rm{fin}}$}
    \STATE \algbreak;
    \ELSE 
    \STATE $\mu=\beta \mu_{k-1}$;
    \STATE $\kappa_{k}=\frac{\kappa_{k-1}}{\beta}$; 
    \WHILE {1}
    \STATE Obtain $\theta_k \in \Omega_{\theta}$ that satisfy $\norm{g_k(\theta_k)-g_k(\theta)} < \mu$;
    \STATE Obtain $u[T_k,T_k+N_k-1]$ that satisfy: 
    \STATE $\bullet$ $N_k \leq N_{x(T_k)}$,
    \STATE $\bullet$ $u[T_k,T_k+N_k-1] \in \overline{B}^{\infty}(0,\rho_{x(T_k)})$,
    \STATE $\bullet$ $\norm{h_k(\theta_k) - x_{\ast}}<\frac{\varepsilon_{\rm{fin}}}{2}$;
    \IF {$h_k(B(\theta_k,\kappa_k \mu)) \subset B(h_k(\theta_k),\frac{\varepsilon_{\rm{fin}}}{2})$}
    \STATE $\mu_k=\mu$; 
    \STATE \algbreak;
    \ELSE
    \STATE $\mu=\beta \mu$;
    \ENDIF
    \ENDWHILE
    \STATE $T_{k+1}=T_{k}+N_k;$
    \STATE Apply $u[T_k,T_k+N_k-1]$ to system (\ref{eq:system}) to obtain $x(T_{k+1})$.  
    \STATE $k=k+1$;
    \ENDIF
  \end{algorithmic}
 \end{alg}
 
 \begin{theorem}
   Under Assumptions~\ref{asm:reachability}, \ref{asm:excitation}, and ~\ref{asm:N_x and rho_x}, 
   Algorithm~\ref{alg:excitation_finite_error}
   terminates after finitely many iterations, 
   and the state of (\ref{eq:system}) reaches to the neighborhood 
   $B(x_{\ast},\varepsilon_{\rm{fin}})$ of the target state $x_{\ast}$.  
 \end{theorem}

 \proof
 We first prove that Algorithm~\ref{alg:excitation_finite_error} is feasible.  
 Assumptions~\ref{asm:reachability} and  \ref{asm:N_x and rho_x} make
 all steps inside the while loop feasible, except for 
 the condition 
 \begin{equation}
   h_k(B(\theta_k,\kappa_k \mu)) \subset B(h_k(\theta_k),\frac{\varepsilon_{\rm{fin}}}{2}).  
   \label{eq:finite_error_termination}
 \end{equation}
 The analysis of $h_k(\theta_k)$ needs some care, 
 because it is the abbreviation of the function 
 $\varphi(T_k+N_k,T_k,x(T_k);u,\theta_k)$.  
 However, since $N_k$ is bounded by $N_{x(T_k)}$, 
 the amplitude of $u$ is bounded by $\rho_{x(T_k)}$, 
 and $\varphi$ is $C^1$, 
 for a positive constant $c(x(N_k))$ that depends on $x(N_k)$, 
 $\norm{h_k(\theta)-h_k(\theta^{\prime})} \leq c(x(N_k)) \norm{\theta-\theta^{\prime}}$
 for all $\theta,\theta^{\prime} \in \Omega_{\theta}$.  
 If (\ref{eq:finite_error_termination}) fails, 
 the minor loop of the while loop of Algorithm~\ref{alg:excitation_finite_error} 
 makes $\mu=\beta \mu$, $0 < \beta < 1$.  
 Thus, residually, $\mu < \frac{\varepsilon_{\rm{fin}}}{2\kappa_k \max\{1,c(x(N_k))\}}$, 
 and (\ref{eq:finite_error_termination}) is fulfilled.  

 Next, we prove by contradiction that Algorithm~\ref{alg:excitation_finite_error} terminates after finitely many steps.  
 Suppose that the termination condition of the (Loop) part of Algorithm~\ref{alg:excitation_finite_error} is never fulfilled.  
 Then, an infinite sequence of the parameter estimate $(\theta_k)_{k \in \Nset}$ is obtained.  
 In this case, the fourth line of (Loop) makes $\mu_k \leq \beta \mu_{k-1}$ (in fact, 
 with the iteration of the while loop, $\mu_k =\beta^{d_{k-1}} \mu_{k-1}$ for some $d_{k-1} \geq 1$. )
 Thus, the sequence $(\mu_k)_{k \in \Nset}$ converges to zero.  
 Contrary, by the execution of the fifth line of (Loop), $\kappa_k=\frac{\kappa_{k-1}}{\beta}$, and $\kappa_0 >0$; 
 hence the sequence $(\kappa_k)_{k \in \Nset}$ diverges to infinity.  
 Because $\theta_k \in  \Omega_{\theta}$ for each $k$ and $\Omega_{\theta}$ is compact, 
 the sequence $(\theta_k)_{k \in \Nset}$ has accumulation points in $\Omega_{\theta}$.  
 Let $L$ be the set of all accumulation points of $(\theta_k)_{k \in \Nset}$.  
 For $q \in L$, there is a subsequence $(\theta_{k_l})_{l \in \Nset}$ that converges to $q$.  
 For all $j$, $\exists k_l \geq j$, and because 
 $g_j$ is an initial segment of $g_{k_l}$ for $j \leq k_l$, 
 \[
 \norm{g_j (\theta_{k_l})-g_j (\theta) )}
 \leq \norm{g_{k_l}(\theta_{k_l})-g_{k_l}(\theta)}
 \leq \mu_{k_l}.  
 \]
 Because $\mu_{k_l}$ converges to zero and $g_j$ is continuous, 
 $\lim_{l\rightarrow \infty}\norm{g_j (\theta_{k_l})-g_j (\theta) )}=\norm{g_j(q)-g_j(\theta)}=0$.  
 Hence, 
 \begin{equation}
   \forall q \in L, \forall j, g_j(q)=g_j(\theta).  
   \label{eq:q_vs_theta}
 \end{equation}

 Next, let $L(\varepsilon)=\cup_{q \in L} B(q,\varepsilon)$ for some $\varepsilon >0$.  
 Then, we can show that 
 \begin{equation}
   \forall \varepsilon, \exists \overline{k}_{a}, \forall k \geq \overline{k}_{a},  \theta_k \in L(\varepsilon).  
   \label{eq:L_neighborhood}
 \end{equation}
 To see this, let us suppose contrary: 
 $\exists \varepsilon, \forall \overline{k}_{a}, \exists k \geq \overline{k}_{a}, \theta_k \not\in  L(\varepsilon)$.  
 Then, $(\theta_k)_{k \in \Nset}$ has an accumulation point in $\Omega_{\theta} \setminus L(\varepsilon)$, 
 contradicting the assumption that $L$ is the set of all accumulation points.  
 Let $\varepsilon_{g_1}$ and $c_{g_1}$ be constants defined in Lemma~\ref{lem:uniform_jacobian}.  
 Choose a $\overline{k}_{a}$ that satisfies (\ref{eq:L_neighborhood}) for $\varepsilon=\varepsilon_{g_1}$.  
 Because $(\kappa_k)_{k \in \Nset}$ diverges to infinity, 
 $\exists \overline{k}_{b}$,  $\forall k \geq \overline{k}_{b}$, 
 $\kappa_k \geq \frac{1}{c_{g_1}}$.  
 Let $k \geq \max\{\overline{k}_a,\overline{k}_b\}$.  
 Let $q_k=\arg \min_{q \in L} \norm{\theta_k-q}$.  
 Then, $\norm{\theta_k-q_k}<\varepsilon_{g_1}$.  
 Hence, by Lemma~\ref{lem:uniform_jacobian}, 
 $\norm{\theta_k-q_k}\leq \frac{1}{c_{g_1}} \norm{g_1(\theta_k)-g_1(q_k)}$, 
 and because  $\kappa_k \geq \frac{1}{c_{g_1}}$, 
 $\norm{\theta_k-q_k}\leq  \kappa_k \norm{g_1(\theta_k)-g_1(q_k)}$.  
 Moreover, by (\ref{eq:q_vs_theta}), $g_1(q_k)=g_1(\theta)$, 
 and with the first execution of the statement inside the while loop of Algorithm~\ref{alg:excitation_finite_error}, 
 $\norm{g_1(\theta_k)-g_1(\theta)} \leq \norm{g_k(\theta_k)-g_k(\theta)} \leq \mu_k$.  
 Thus,  $\norm{\theta_k-q_k}\leq  \kappa_k \mu_k$, 
 and hence,
 $h_k(q_k) \subset B(h_k(\theta_k),\frac{\varepsilon_{\rm{fin}}}{2})$.  
 Because $h_k(\theta_k) \in B(x_{\ast},\frac{\varepsilon_{\rm{fin}}}{2})$, 
 $h_k(q_k) \in B(x_{\ast},\varepsilon_{\rm{fin}})$, 
 which shows that the termination condition of (Loop) has already been fulfilled  at the $k+1$-th step; a contradiction.  
 \qed

 \section{Conclusion}
  In this paper, 
  the finite-time adaptive set-point regulation problem 
  for discrete-time nonlinear systems with parametric uncertainty
  has been solved under the assumption that the target state is reachable provided 
  that the parameter is known and an excitation signal is available.  

  The proposed controller has a pathological structure that
  all history of the past state is preserved until the 
  state reaches to the target state.  
  Moreover, in Algorithm\ref{alg:excitation_finite_error}, it is not easy to
  obtain estimates of $N_{x(T_k)}$ and $\rho_{x(T_k)}$.  
  Hence, the proposed algorithms are computationally extremely demanding 
  and by no means practical.  
  They should be regarded as being of purely theoretical and conceptual nature.  
  On the other hand, 
  in order for the proposed algorithms to be applicable, 
  except for the reachability to the target state,
  no additional condition is required. 
  This contrasts to the majority of existing methods of nonlinear adaptive control, 
  where many structural conditions are required in order for those methods to be applicable.  
  Thus, the implication of this paper is to show the potential of the concept of nonlinear adaptive control
  in the sense that 
  no extra condition other than the reachability to the target state is required in order to construct a stable nonlinear adaptive controller
  provided that sufficiently fast and reliable nonlinear minimizer or nonlinear equation solver is available.  
  In this respect, nonlinear adaptive control problem is reduced to nonlinear optimization problem.  
  An extra bonus of the proposed methodology is that it is completely ``Lyapunov-free.''  

  From a practical point of view, it is desirable to develop a more down-to-earth algorithm that has less computational complexity
  but does not necessitate extra conditions other than the reachability to the target state.  
  It is also to be noted that the proposed algorithms have the drawback that they are not robust against disturbances.  
  To overcoming these problems is left for further research.  


\end{document}